\documentclass[11pt]{article}

\usepackage[english]{babel}
\usepackage{geometry}
\usepackage[utf8]{inputenc}
\usepackage{mathtools,amsfonts,amsthm,mathrsfs,amssymb,stmaryrd}
\usepackage{textcomp}
\usepackage{bookmark}
\usepackage{microtype}
\usepackage{enumitem}
\usepackage{todonotes}
\usepackage{fullpage}
\usepackage{verbatim}
\usepackage{cleveref}
\usepackage{authblk}
\usepackage{thm-restate}
\usepackage{xstring,refcount}

\usepackage{tikz}
\usetikzlibrary{arrows.meta, shapes.misc, positioning}

\newtheorem{theorem}{Theorem}[section]
\newtheorem{lemma}[theorem]{Lemma}

\newtheorem{proposition}[theorem]{Proposition}

\newtheorem{corollary}[theorem]{Corollary}
\newtheorem{conjecture}{Conjecture}
\newtheorem{question}{Question}

\newtheorem{claim}{Claim}[theorem]

\theoremstyle{definition}
\newtheorem*{definition*}{Definition}

\newenvironment{proofclaim}[1][]{\par\noindent {\it Proof of claim}. }{ \hfill$\lozenge$\par\addvspace{6pt plus 6pt}}

\newcommand{\bI}{\mathbf{I}}

\newcommand{\discr}{\varphi}

\newcommand{\C}{\mathscr{C}}

\newcommand{\pr}[1]{\mathbb{P}\left[#1\right]}

\newcommand{\sset}[2]{\left\{#1 : #2 \right \}}
\newcommand{\pth}[1]{\left(#1\right )}

\newcommand{\floor}[1]{\left\lfloor #1 \right\rfloor}
\newcommand{\ceil}[1]{\left\lceil #1 \right\rceil}
\newcommand{\induced}{\subseteq_{\rm I}}
\newcommand{\Wlog}{Without loss of generality}
\renewcommand{\emptyset}{\varnothing}

\newcommand{\Myc}[2]{\mathscr{M}_{#1}^{#2}}
\newcommand{\dijoin}{\wedge}

\author[1]{Timoth\'ee Corsini}
\author[2]{Lucas Picasarri-Arrieta\footnote{Research supported by JST as part of ASPIRE, Grant Number JPMJAP2302.}$^,$}
\author[3]{Th\'eo Pierron}
\author[4]{Fran\c{c}ois Pirot}
\author[5]{Eileen Robinson}

\affil[1]{Université de Bordeaux, LaBRI, France.}
\affil[2]{National Institute of Informatics, Tokyo, Japan.}
\affil[3]{Université Lyon 1, LIRIS UMR CNRS 5205, F-69621, Lyon, France.}
\affil[4]{Université Paris-Saclay, France.}
\affil[5]{Université libre de Bruxelles, Belgium.}

\title{Chromatic discrepancy of locally $s$-colourable graphs}
\date{}

\begin{document}

\maketitle

\begin{abstract}
    The chromatic discrepancy of a graph $G$, denoted $\discr(G)$, is the least over all proper colourings $\sigma$ of $G$ of the greatest difference between the number of colours $|\sigma(V(H))|$ spanned by an induced subgraph $H$ of $G$ and its chromatic number $\chi(H)$. We prove that the chromatic discrepancy of a triangle-free graph $G$ is at least $\chi(G)-2$. This is best possible and positively answers a question raised by Aravind, Kalyanasundaram, Sandeep, and Sivadasan.

    More generally, we say that a graph $G$ is locally $s$-colourable if the closed neighbourhood of any vertex $v\in V(G)$ is properly $s$-colourable --- in particular, a triangle-free graph is locally $2$-colourable. We conjecture that every locally $s$-colourable graph $G$ satisfies $\discr(G) \geq \chi(G)-s$, and show that this would be almost best possible.
    We prove the conjecture when $\chi(G) \le 11s/6$, and as a partial result towards the general case, we prove that every locally $s$-colourable graph $G$ satisfies $\discr(G) \geq \chi(G) - s\ln \chi(G)$. 
    If the conjecture holds, it implies in particular, for every integer $\ell\geq 2$, that any graph $G$ without any copy of $C_{\ell+1}$, the cycle of length $\ell+1$, satisfies $\discr(G) \geq \chi(G) - \ell$. When $\ell \ge 3$ and $G\neq K_\ell$, we conjecture that we actually have $\discr(G)\ge \chi(G) - \ell + 1$, and prove it in the special case $\ell = 3$ or $\chi(G) \le 5\ell/3$. In general, we further obtain that every $C_{\ell+1}$-free graph $G$ satisfies $\discr(G) \geq \chi(G) - O_{\ell}(\ln \ln \chi(G))$. We do so by determining an almost tight bound on the chromatic number of balls of radius at most $\floor{\ell/2}$ in $G$, which could be of independent interest. 
\end{abstract}

\section{Introduction}

A proper colouring of a (simple, undirected) graph $G$ is an assignment of colours to the vertices of $G$ that induces no monochromatic edge. 
Understanding which graphs need many colours in a proper colouring is a longstanding open question. In~\cite{aravind2015}, the authors introduce \emph{chromatic discrepancy}, a notion that captures how far from optimal a proper colouring of $G$ can be when restricted to its induced subgraphs. More precisely, given a non-empty graph $G$ and a proper colouring $\sigma$ of $G$, the \emph{discrepancy of $\sigma$} is
\[ \discr_{\sigma}(G)\coloneqq \max_{H\induced G} \;   \Big(|\sigma(V(H))| - \chi(H) \Big),\]
where $H\induced G$ stands for $H$ being an induced subgraph of $G$, and $\chi(H)$ denotes the chromatic number of $H$. The \emph{chromatic discrepancy of $G$} is the minimum discrepancy over all proper colourings of $G$, denoted
\[ \discr(G)\coloneqq \min_{\sigma \in \C(G)} \discr_{\sigma}(G),\]
where $\C(G)$ denotes the set of all proper colourings of $G$.
Note that, for any graph $G$ and any proper $\chi(G)$-colouring $\sigma$ of $G$, we have $|\sigma(V(H))|-\chi(H) \le \chi(G)-1$ for every $H\induced G$, that is $\discr_\sigma(G)\le \chi(G)-1$. This shows that every graph $G$ satisfies $\discr(G) \le \chi(G)-1$, and for equality to hold one must in particular be able to find an independent set that spans every colour in every optimal proper colouring of $G$.
However, we note that this condition is not sufficient in general, as it can be seen, for instance, with the construction in \cite[Theorem 7.1]{aravind2015} which essentially consists of the Mycielskian of a complete graph, minus the universal vertex for the twins. In such a graph, the discrepancy of an optimal colouring is $1$, while the discrepancy of the graph is arbitrarily large.

The goal of this paper is to provide a collection of new lower bounds on $\discr(G)$. Observe that proving a bound of the form $\discr(G)\ge k$ amounts to finding, for every proper colouring of $G$, an induced subgraph $H$ of $G$ receiving at least $\chi(H)+k$ colours, that is $k$ more colours than needed to colour $H$. In~\cite{aravind2015}, the authors raised the following question that links the chromatic discrepancy of a graph $G$ and its clique number $\omega(G)$.

\begin{question}
\label{question:aravind}
    Does $\discr(G)\ge \chi(G)-\omega(G)$ hold for every graph $G$?
\end{question}

As a partial answer to this question, they obtained that every graph $G$ has chromatic discrepancy at least $\discr(G) \geq \frac{1}{2} \left( \chi(G) - \omega(G)\right)$, and showed that if true, the bound of Question~\ref{question:aravind} would be best possible when $\omega(G)=2$, that is when $G$ is a triangle-free graph.

A negative answer to Question~\ref{question:aravind} was provided in~\cite[Corollary~4]{aravind2021structure}. It is even proved that, for every $k \ge 3$, there exists a $K_4$-free graph $G$ with chromatic number $\chi \geq k$ and
\begin{equation}
    \label{eq:gap_fixed_omega}
    \discr(G)\leq \chi-\Omega\pth{\frac{\chi^{1/3}}{(\log \chi)^{4/3}}}.
\end{equation}

This leaves the triangle-free case open for upper bound. They provide the lower bound $\discr(G)\ge \chi(G)-\log_2 \chi(G)-1$ for this case, and ask whether the logarithmic term could be replaced by a constant. Our first main contribution is a positive answer to this question.

\begin{theorem}
\label{thm:main}
    Every triangle-free graph $G$ satisfies $\discr(G)\ge \chi(G)-2$. 
\end{theorem}

Note that, since every graph $G$ has chromatic discrepancy at most $\chi(G)-1$, it follows from Theorem~\ref{thm:main} that every triangle-free graph $G$ has chromatic discrepancy either $\chi(G)-2$ or $\chi(G)-1$.

In~\cite{aravind2015}, the chromatic discrepancy is also linked to the notion of \emph{local colouring}, whose purpose is to minimise the number of colours used on the neighbourhood of each vertex in a proper colouring. 
The \emph{local chromatic number} of $G$ is
\[
    \psi(G) \coloneqq \min_{\sigma\in \C(G)} \max_{v\in V(G)} |\sigma(N[v])|. 
\]

Given integers $r,s\ge 1$, a graph $G$ is \emph{$r$-locally $s$-colourable} if $\chi(B_r(v))\le s$ for every $v\in V(G)$, where $B_r(v)$ is the \emph{ball of radius $r$ centred in $v$ in $G$}, that is, the subgraph of $G$ induced by the vertices at distance at most $r$ from $v$. When $r$ is omitted, its implicit value is $1$.
Note that, if a graph $G$ has local chromatic number $\psi(G)\leq s$, then it is locally $s$-colourable, but the converse is not true in general.

It easily follows from the definitions that, for every integer $s\ge 2$ and every locally $s$-colourable graph $G$, we have
\begin{equation}
    \label{eq:lowerbound_psi}
    \discr(G) \geq \psi(G) - s.
\end{equation}

A nice consequence of~\eqref{eq:lowerbound_psi} is that the chromatic discrepancy of a locally $s$-colourable graph $G$ is related to its fractional chromatic number $\chi_f(G)$, namely 
\begin{equation}
    \label{eq:chif-discr}
    \discr(G) \ge \chi_f(G)-s.
\end{equation}

Let us recall that the \emph{fractional chromatic number $\chi_f(G)$} of a given graph $G$ is the minimum $w$ such that there is a probability distribution over the independent sets of $G$ such that, for each vertex $v$, given an independent set $\bI$ drawn from that distribution, $\bI$ contains $v$ with probability at least $1/w$.
We note that there are many other equivalent definitions for the fractional chromatic number of a graph, the aforementioned one being the most convenient for our purpose. 
Equation~\eqref{eq:chif-discr} follows from~\eqref{eq:lowerbound_psi} and the following classical result, of which we repeat the short proof for completeness. 

\begin{lemma}
    \label{lem:bounded_fractional_chromatic_number}
    For every graph $G$, $\chi_f(G) \le \psi(G)$.
\end{lemma}

\begin{proof}
    Let $\sigma$ be a proper colouring of $G$ such that $|\sigma(N[v])| \le k$ for every vertex $v\in V(G)$. 
    Let $\prec$ be a uniformly random ordering of the colours of $\sigma$. We construct a random independent set $\bI$ by adding a vertex $v$ in $\bI$ if $\sigma(v) \prec \sigma(u)$ for every $u\in N(v)$.
    It is straightforward that $\pr{v\in \bI}\ge 1/k$ for every $v\in V(G)$, and so $\chi_f(G)\le k$ by definition. 
\end{proof}

In contrast, it is known that there exist graphs with arbitrarily large chromatic number and with local chromatic number at most $3$; shift graphs~\cite{erdosTG68} (defined on triplets) are such graphs. 
We pose the following conjecture that states a non-fractional version of \eqref{eq:chif-discr}, and is half-way between Question~\ref{question:aravind} and Equation~\eqref{eq:lowerbound_psi} (since $\chi(G) \geq \psi(G)$ and $s\geq \omega(G)$). 

\begin{conjecture}
    \label{conjecture:main_conj}
    For every integer $s\ge 2$, every locally $s$-colourable graph $G$ satisfies
    \[\discr(G)\ge \chi(G)-s.\]
\end{conjecture}

We note that, if true, \Cref{conjecture:main_conj} is almost best possible, by constructing, for all integers $2\le s\le k$, a locally $s$-chromatic graph $G$ of chromatic number $k$  with $\discr(G) = k-s+1$ (see \Cref{prop:no_rainbow_IS_mycielski}).

Observe that \Cref{thm:main} corresponds to the case $s=2$ of \Cref{conjecture:main_conj}, since the locally $2$-colourable graphs are exactly the triangle-free graphs. 
In \cite{aravind2015}, the authors characterised all graphs $G$ with $\discr(G)=0$ as complete multipartite graphs. This implies that \Cref{conjecture:main_conj} holds when $\chi(G)=s+1$. 
We extend this result by showing that \Cref{conjecture:main_conj} holds whenever $\chi$ is sufficiently close to $s$.

\begin{theorem}
\label{thm:conj-special-case}
    For every integer $s\ge 2$, every locally $s$-colourable graph $G$ of chromatic number $\chi < \frac{11s+16}{6}$ satisfies
    \[\discr(G)\ge \chi-s.\]
\end{theorem}

For general values of $s$ and $\chi$,  a generalisation of the arguments used in \cite{aravind2021structure} yields the following partial results towards \Cref{conjecture:main_conj}.

\begin{theorem}
    \label{thm:locally-s-chromatic}
    For every integer $s\ge 3$, every locally $s$-colourable graph $G$ satisfies
    \[\discr(G)\ge \chi(G)-s \ln \chi(G).\]
\end{theorem}

We note that no analogue of Theorem~\ref{thm:locally-s-chromatic} holds when $s$ is replaced by $\omega(G)$. Indeed, recall that there exists a $K_4$-free graph $G$ with arbitrarily large chromatic number $\chi$ satisfying~\eqref{eq:gap_fixed_omega}. In particular, if $\chi$ is chosen large enough, then such a graph $G$ satisfies
\[
    \discr(G) \leq \chi - f(\omega) \ln (\chi)
\]
for any fixed function $f\colon \mathbb{N}\to \mathbb{N}$. 

In the more restricted case of $2$-locally $s$-colourable graphs, we obtain the following stronger bound.

\begin{theorem}
    \label{thm:2-local}
    For all integers $s_1,s_2 \ge 2$, there exists an absolute constant $C\ge 0$ such that the following holds.
    Every $1$-locally $s_1$-colourable and $2$-locally $s_2$-colourable graph $G$ satisfies
    \[
    \discr(G) \ge \chi(G) - 4s_1s_2\ln \ln \chi(G) - C.
    \]
\end{theorem}

We finally propose a weaker version of Conjecture~\ref{conjecture:main_conj}. Given an integer $\ell \geq 2$, a $C_{\ell+1}$-free graph is any graph that contains no cycle of length exactly $\ell+1$ as a (not necessarily induced) subgraph. It is straightforward, see Proposition~\ref{prop:C_ell_Nv_degenerate}, that every $C_{\ell+1}$-free graph is locally $\ell$-colourable. The following weaker form of Conjecture~\ref{conjecture:main_conj} is open.

\begin{conjecture}
    \label{conjecture:weaker_conj}
    For every integer $\ell \ge 2$, every $C_{\ell+1}$-free graph $G$ satisfies 
    \[ 
        \discr(G)\ge \chi(G)-\ell.
    \] 
\end{conjecture}

We suspect that equality may occur only when $\ell=2$ or when $G$ is a complete graph, which leads us to propose the following stronger conjecture.

\begin{conjecture}
    \label{conj:weak_conj}
    For every integer $\ell \ge 3$, every $C_{\ell+1}$-free graph $G \neq K_{\ell}$ satisfies 
    \[ 
        \discr(G)\ge \chi(G)-\ell+1.
    \] 
\end{conjecture}

Recall that Conjecture~\ref{conjecture:weaker_conj} holds when $\ell=2$, as it corresponds to Theorem~\ref{thm:main}.
As evidences for Conjecture~\ref{conj:weak_conj}, we prove that it holds when $\ell=3$ or when $\chi(G)$ is sufficiently close to $\ell$.

\begin{theorem}
    \label{thm:C4_free}
    Every $C_4$-free graph $G\neq K_3$ satisfies
    \[\discr(G)\ge \chi(G)-2.\]
\end{theorem}

\begin{theorem}
    \label{thm:weak_conj_special_case}
    For every integer $\ell \ge 3$, every $C_{\ell+1}$-free graph $G\neq K_\ell$ of chromatic number $\chi < \frac{5\ell+2}{3}$ satisfies
    \[ 
    \discr(G) \ge \chi - \ell + 1. 
    \]
\end{theorem}

One of the main differences between $C_{\ell+1}$-free graphs and locally $\ell$-colourable graphs is that forbidding $C_{\ell+1}$ yields some information at distance more than $1$ from each vertex. In particular, as we show in Section~\ref{sec:coloring-balls}, the balls of radius $\floor{\ell/2}$ are $2\ell$-colourable (see Theorem~\ref{thm:cl}). 
Combined with Theorem~\ref{thm:2-local}, this implies the following partial result towards the general statements of Conjectures~\ref{conjecture:weaker_conj} and~\ref{conj:weak_conj}.

\begin{theorem}
    \label{thm:discr_c_ell}
    For every integer $\ell \ge 4$, every $C_{\ell+1}$-free graph $G$ satisfies
    \[ \discr(G) \ge \chi(G) - O_{\ell}(\ln \ln \chi(G)).\]
\end{theorem}

\paragraph{Outline of the paper.} We first recall some classical notation and well-known results, and introduce specific definitions in Section~\ref{sec:preliminaries}.
In \Cref{sec:triangle-free}, we prove \Cref{thm:main} and show that it is tight. 
In \Cref{sec:discr_locally_s_colourable}, we prove that \Cref{conjecture:main_conj}, if true, is almost best possible. We then prove Theorem~\ref{thm:conj-special-case} and provide a short proof of Theorem~\ref{thm:locally-s-chromatic}.
In \Cref{sec:2-local}, we prove that the chromatic discrepancy of $2$-locally $s$-colourable graphs $G$ is $\chi(G) - O_s(\ln\ln\chi(G))$, thus deriving Theorem~\ref{thm:2-local}.
In \Cref{sec:cl-free}, which is dedicated to graphs excluding a cycle, we prove Theorems~\ref{thm:C4_free} and~\ref{thm:weak_conj_special_case}. 
In \Cref{sec:coloring-balls}, we show that balls of small radius have bounded chromatic number in $C_{\ell+1}$-free graphs, and finally derive a proof of \Cref{thm:discr_c_ell}.

\section{Preliminaries}
\label{sec:preliminaries}

\subsection{Terminology and notation}

Given a vertex $v$ of a graph $G$ and an integer $r \ge 1$, we define $N^r(v)$ (respectively $N^r[v]$) as the set of vertices at distance exactly $r$ (respectively at most $r$) from $v$ in $G$. Also, we define the \emph{layer at distance $r$ from $v$}, denoted by $L_r(v)$, as the subgraph of $G$ induced by $N^r(v)$. Analogously, the \emph{ball of radius $r$ with centre $v$}, denoted by $B_r(v)$, is the subgraph of $G$ induced by $N^r[v]$. 

Given an integer $k$, we say that $G$ is \emph{$k$-degenerate} if every non-empty subgraph $H$ of $G$ contains a vertex $v$ of degree at most $k$ in $H$. We will use many times the easy observation that a $k$-degenerate graph has chromatic number at most $k+1$. In particular, a graph with chromatic number $\chi$ contains an induced subgraph with minimum degree at least $\chi-1$.

For every integer $\ell$, we denote by $P_\ell$ the path on $\ell$ vertices. Further, for any $\ell\geq 3$, we denote by $C_\ell$ the cycle on $\ell$ vertices. The {\it length} of a path or cycle is its number of edges.

For every integer $n\ge 1$, we let $[n]$ denote the set of integers $\{1, \ldots, n\}$.

\subsection{Basic properties and well-known results}

In this subsection, we provide basic properties and well-known results used later on.
We first justify that $\discr$ is non-increasing under taking induced subgraphs.

\begin{proposition}
    Let $G$ be a graph and $H$ be an induced subgraph of $G$, then $\discr(G) \geq \discr(H)$. 
\end{proposition}
\begin{proof}
    Let $\sigma$ be a proper colouring of $G$ such that $\discr(G) = \discr_{\sigma}(G)$. Let $\sigma_H$ be the restriction of $\sigma$ to $V(H)$. Since every induced subgraph of $H$ is an induced subgraph of $G$, it follows from the definitions that
    \[
        \discr(H) \leq \discr_{\sigma_H}(H) \leq \discr_{\sigma}(G) = \discr(G). \qedhere
    \]
\end{proof}

We make use of the following property of $C_{\ell+1}$-free graphs many times.

\begin{proposition}
    \label{prop:C_ell_Nv_degenerate}
    Let $G$ be a $C_{\ell+1}$ free graph and let $v$ be any vertex of $G$. Then $G[N(v)]$ is $(\ell-2)$-degenerate.
\end{proposition}
\begin{proof}
    Assume that this is not the case, so $G[N(v)]$ contains a subgraph $H$ with minimum degree $\ell-1$. Then any maximum path in $H$ has length at least $\ell$, which together with $v$ yields a copy of $C_{\ell+1}$ in $G$, a contradiction.
\end{proof}

We finally make use of the so-called Brooks' Theorem.

\begin{theorem}[Brooks~\cite{brooksMPCPS37}]
    \label{thm:brooks}
    Let $G$ be a connected graph of maximum degree $\Delta$ and chromatic number $\Delta+1$. Then either $\Delta=2$ and $G$ is an odd cycle, or $G = K_{\Delta+1}$.
\end{theorem}

\subsection{A variant of the chromatic discrepancy}

Our proof of Theorem~\ref{thm:main} actually relies on a graph parameter whose definition is quite similar to that of chromatic discrepancy. 
The first motivation for introducing our parameter follows from the observation that most of the proofs from~\cite{aravind2015,aravind2021structure} rely on finding some so-called \emph{rainbow} structures and using them to bound the chromatic discrepancy. With this motivation in mind, we define the following. Given a graph $G$ and an integer $\chi(G) \le p \le |V(G)|$, we let $\C_p(G)$ denote the set of proper colourings of $G$ that use exactly $p$ colours --- we note that this set is non-empty for every $\chi(G) \le p \le |V(G)|$, and would be empty for other values of $p$.
Given a proper $p$-colouring $\sigma \in \C_p(G)$, a subset of vertices $X\subseteq V(G)$ is said to be \emph{$\sigma$-rainbow} if the colours $\sigma(x)$ are distinct over all $x\in X$ (so $|\sigma(X)|=|X|$). A $\sigma$-rainbow set of size $p$ is called a \emph{rainbow cover of $\sigma$}.
We are interested in the minimum chromatic number of a rainbow cover of $\sigma$. To that end, 
for every $\chi(G) \le p \le |V(G)|$, we define

\begin{equation}
    \label{eq:def-f_G}
f_G(p) \coloneqq \max_{\sigma \in \C_p(G)} \min \big\{ \chi(G[X]) : X \subseteq V(G), |\sigma(X)|=p\big\}.
\end{equation}

A second motivation for introducing this parameter comes from the observation that, in the definition of chromatic discrepancy, considering colourings $\sigma$ with at least $2\chi(G)$ colours is useless, since then \[\max\limits_{H\subseteq_I G} |\sigma(V(H))|-\chi(H) \ge |\sigma(V(G))|-\chi(G) \geq \chi(G) > \discr(G).\]
In contrast, our parameter is designed to take into account all possible proper colourings of $G$ using $p$ colours, regardless of their optimality (that is, $p$ can be arbitrarily far from $\chi(G)$). 
Our parameter is related to the chromatic discrepancy as follows.

\begin{lemma}
\label{lem:discrfp}
    For every graph $G$, 
    \[\discr(G)=\min  \big\{p-f_G(p) : \chi(G)\le p\le |V(G)|\big\}.\]
\end{lemma}

\begin{proof}
    Let $\sigma \in \C(G)$ and $H\subseteq_I G$ be such that $\discr(G) = |\sigma(V(H))| - \chi(H)$, and let $p$ be the number of colours used by $\sigma$. We first prove that $\discr(G) \geq p-f_G(p)$, which implies $\discr(G) \ge \min_p (p-f_G(p))$.
    
    By definition of $f_G(p)$, $G$ contains an induced subgraph $H'$ such that $V(H')$ is a rainbow cover of $\sigma$ and $\chi(H') \leq f_G(p)$. As desired, by choice of $H$, we obtain
    \[
    \discr(G) = |\sigma(V(H))| - \chi(H) \geq |\sigma(V(H'))| - \chi(H') \geq p - f_G(p). 
    \]

    To complete the proof of the statement, it is sufficient to show that $\discr(G) \leq p-f_G(p)$ holds for every integer $p\in \{\chi(G),\ldots, |V(G)|\}$. Let us thus fix $p$, and let $\sigma \in \C_p(G)$ be such that $\min \{ \chi(G[X]) : |\sigma(X)|=p\}$ is maximised, so that $f_G(p) = \min \{ \chi(G[X]) : |\sigma(X)|=p\}$. Let $H \subseteq_I G$ be such that $|\sigma(V(H))| - \chi(H)$ is maximised. 
    We may assume that $H$ spans all $p$ colours of $\sigma$, for otherwise we may add to $H$ one vertex of each missing colour in $\sigma$, and obtain a graph $H'$ spanning all colours of $\sigma$ such that $|\sigma(V(H'))| - \chi(H') \geq |\sigma(V(H))| - \chi(H)$.

    By definition of $\discr(G)$, we have $\discr(G) \leq |\sigma(V(H))| - \chi(H) = p - \chi(H)$. On the other hand, $\chi(H) \geq \min \{\chi(H') : |\sigma(H')|=p\}=f_G(p)$.
    We deduce $\discr(G) \leq p -f_G(p)$ as desired.
\end{proof}

As a consequence of Lemma~\ref{lem:discrfp}, proving, for every integer $\chi(G) \leq p \leq |V(G)|$, an upper bound of the form $f_G(p)\le p-C$ for some constant $C>0$ yields the lower bound $\discr(G)\ge C$. This is how we prove Theorem~\ref{thm:main} in the next section. 

\section{Chromatic discrepancy of triangle-free graphs}
\label{sec:triangle-free}

The goal of this section is to prove Theorem~\ref{thm:main}, that we first recall here for convenience.
\begin{theorem}
    Every triangle-free graph $G$ satisfies $\discr(G)\ge \chi(G)-2$. 
\end{theorem}

This theorem is claimed to be tight in~\cite{aravind2015} for the Mycielski graphs. For completeness, we first include a proof of this result. 

\subsection{A family of triangle-free graphs \texorpdfstring{$G$}{G} with chromatic discrepancy \texorpdfstring{$\chi(G)-2$}{χ(G)-2}}

Recall that the Mycielski graphs are defined by taking $M_2$ as the 1-edge graph, and letting $M_k$ be the \emph{Mycielskian} of $M_{k-1}$ for every $k\ge 3$, that is $M_k$ is obtained from $M_{k-1}$ by adding a twin of each vertex and a vertex $v_0$ adjacent to all the twins. It is well-known that each graph $M_k$ is triangle-free and has chromatic number $k$~\cite{mycielski1955}. For every $k\geq 2$, $M_k$ comes with a \emph{canonical} proper $k$-colouring built by colouring $v_0$ with $k$, and each vertex of $M_{k-1}$ and its twin with its canonical colour in the canonical $(k-1)$-colouring of $M_{k-1}$.

\begin{proposition}
    \label{prop:no_rainbow_IS_mycielski}
    Let $G=M_k$ be the $k$-th iteration of the Mycielski construction, and let $\sigma\in \C_k(G)$ be its canonical colouring. Then there is no independent set spanning all $k$ colours in $G$.
\end{proposition}

\begin{proof}
    We prove the result by induction on $k$. The first graph of the construction is $G=K_2$, where every independent set has size $1$, so the initialisation of the induction holds.
    Let us now assume that $k\ge 3$. Assume for the sake of contradiction that there is an independent set $I$ of $M_k$ that spans all $k$ colours of its canonical colouring $\sigma$. Then we have $v_0\in I$, and so $I$ contains no twin vertex. Hence $I\setminus\{ v_0\} \subseteq V(M_{k-1})$; this is an independent set that spans $k-1$ colours in the canonical colouring of $M_{k-1}$, a contradiction of the induction hypothesis. This ends the proof.
\end{proof}

It follows from Proposition~\ref{prop:no_rainbow_IS_mycielski} that, in the canonical colouring of $M_k$, every subgraph $H$ of $M_k$ either has $\chi(H)>1$ or spans at most $k-1$ colours. Therefore $\discr(M_k)\le k-2=\chi(M_k)-2$. This shows the tightness of Theorem~\ref{thm:main}.

\subsection{Triangle-free graphs \texorpdfstring{$G$}{G} have chromatic discrepancy at least \texorpdfstring{$\chi(G)-2$}{χ(G)-2}}

The proof of Theorem~\ref{thm:main} is based on the following result, where $f_G$ is the function defined in \eqref{eq:def-f_G}, i.e. $f_G(p)$ is the maximum over all proper $p$-colourings $\sigma$ of $G$ of the minimum chromatic number of a rainbow cover of $\sigma$, that is an induced subgraph $H \subseteq_I G$ such that $|\sigma(V(H))|=p$.

\begin{theorem}\label{thm:chi_plus_k}
    Let $s\ge 2$, $k \ge 0$ be two integers. Every locally $s$-colourable graph $G$ satisfies 
    \[f_G(\chi (G) + k) \le (s-1)(k+1) + 1.\]
\end{theorem}

Observe that, for a triangle-free graph $G$, one can take $s=2$ and obtain from Theorem~\ref{thm:chi_plus_k} that $f_G(\chi(G)+k)\le k+2$ for every fixed $k$. In other words, we get $f_G(p)\le p-\chi(G)+2$ for every fixed~$p$. Hence, by Lemma~\ref{lem:discrfp}, 
\[
\discr(G) = \min_p (p-f_G(p))\ge \chi(G)-2,
\]
and Theorem~\ref{thm:main} follows.

The proof of Theorem~\ref{thm:chi_plus_k} heavily relies on the following lemma. 

\begin{lemma}
    \label{lemma:rainbow_closed_neighbourhood}
    Let $G=(V,E)$ be a graph, let $p =\chi(G)+k$ be an integer with $k\geq 0$, and let $\sigma \in \C_p(G)$ be a proper $p$-colouring of $G$.
    Then there exists $X \subseteq V(G)$ such that $|X| \le k+1$ and $|\sigma(N[X])| = p$. 
    Moreover, for every colour class $U$ of $\sigma$, there exists such a set $X$ such that $X\cap U \neq \emptyset$.
\end{lemma}

\begin{proof}
    We denote by $[p]$ the colours used by $\sigma$ and let $i\in [p]$ be such that $U=\sigma^{-1}(i)$.

    We proceed by induction on $k$. For the sake of better readability, we denote $\chi(G)$ by $\chi$.
    When $k=0$, $\sigma$ is a proper $\chi$-colouring of $G$. If each vertex $v$ coloured with $i$ misses a colour from $[\chi] \setminus \{i\}$ in its neighbourhood, then we can recolour all of them and obtain a proper $(\chi-1)$-colouring of $G$, a contradiction. Therefore, there exists a vertex $v$ coloured $i$ such that $\sigma(N[v]) = [\chi]$ and we can take $X=\{v\}$.

    Assume now that $k\ge 1$, and let $\sigma\colon V(G) \xrightarrow{} [\chi +k]$ be any proper $(\chi +k)$-colouring of $G$. For each colour $j\in [\chi+k]$, we denote by $V_j$ the set of vertices coloured $j$ in $\sigma$. Free to relabel the colours, we may assume that $i=\chi+k$, that is $U=V_{\chi+k}$. If some vertex $v\in V_{\chi+k}$ satisfies $|\sigma(N(v))| \ge \chi+k-1$, then we can set $X=\{v\}$ and we are done.

    We thus assume that every vertex $v\in V_{\chi + k}$ satisfies $|\sigma(N(v))| < \chi +k - 1$. In particular, in $\sigma$, every vertex coloured $\chi+k$ can be recoloured to a colour of $[\chi+k-1]$. Let $\sigma'$ be the proper $(\chi+k-1)$-colouring of $G$ obtained from $\sigma$, where we recolour each vertex of $V_{\chi+k}$ with the smallest available colour. Formally,
    \[\sigma'(v) =\begin{cases}
         \sigma(v) &\text{ if } v\notin V_{\chi+k}\text{,} \\
         \min \;\{ j\in [\chi +k-1] : j\notin \sigma(N(v)) \} &\text{ otherwise.}
    \end{cases}\]

    Let $j$ be the largest colour of $\sigma'(V_{\chi+k})$. By choice of $\sigma'$, observe that every vertex $v \in V_{\chi+k}$ coloured $j$ in $\sigma'$ satisfies $[j-1] \subseteq \sigma(N(v))$.

    By induction, with $\sigma'$ and $j$ playing the roles of $\sigma$ and $i$ respectively, there exists $X' \subseteq V$ such that $|X'| \le k$, $j\in \sigma'(X')$ and $\sigma'(N[X']) = [\chi + k -1]$.

    Let $u \in X'$ be a vertex coloured $j$ in $\sigma'$. If $u\in V_{\chi+k}$, let $v$ be any vertex of $V_j$. Otherwise, we know that $u\in V_j$, and we let $v$ be any vertex of $V_{\chi+k}$ such that $\sigma'(v) = j$. 
    In both cases, we set $X = X' \cup \{v\}$. 
    Note that we indeed have $|X|= |X'|+1\le k+1$ and $\chi+k \in \sigma(X)$ (since $\chi+k\in \sigma(\{u,v\})$). It remains to show that $N[X]$ intersects every colour class of $\sigma$. Observe first that $\sigma(\{u,v\}) = \{j, \chi + k\}$, so $N[X]$ intersects both colours $j$ and $\chi + k$. Moreover, by definition of $\sigma'$, either $u$ or $v$ cannot be recoloured with a colour smaller than $j$, hence $N(\{u,v\})$ intersects every such colour class. Finally, for every colour $z \in [j+1,\chi+k-1]$, every vertex coloured with $z$ in $\sigma'$ is also coloured $z$ in $\sigma$. In particular, $N[X']$ contains a vertex coloured with $z$ in $\sigma$ by induction, and so does $N[X]$, which concludes the proof of the lemma.
\end{proof}

Observe that Lemma~\ref{lemma:rainbow_closed_neighbourhood} is best possible. Indeed, for all integers $\chi,k\ge 0$, there exists a graph $G$ with chromatic number $\chi$ and a proper $(\chi + k)$-colouring $\sigma$ of $G$ such that every set of $k$ vertices misses at least one colour of $\sigma$ in its closed neighbourhood. To see this, let $G'$ be any graph with chromatic number $\chi$. Let $G$ be the graph obtained from $G'$ by adding an independent set $Y=\{u_1,\dots,u_k\}$ of size $k$. Let $\sigma' : V(G') \xrightarrow{} [\chi]$ be an optimal proper colouring of $G'$, and let $\sigma$ be the extension of $\sigma'$ to $G$ obtained by setting $\sigma(u_i) = \chi + i$. Then for every set $X$ of size $k$, either $X \not\subseteq Y$ and $N[X]$ misses a colour from $\sigma(Y)$, or $X=Y$ and $N[X]$ misses colours $[\chi]$.

\medskip

We are now ready to prove Theorem~\ref{thm:chi_plus_k}.
\begin{proof}[Proof of Theorem~\ref{thm:chi_plus_k}]
    This is a consequence of Lemma~\ref{lemma:rainbow_closed_neighbourhood}. Let us fix $\sigma$ a proper $(\chi(G)+k)$-colouring of $G$. By Lemma~\ref{lemma:rainbow_closed_neighbourhood}, let $X\subseteq V(G)$ such that $|X| \le k+1$, and $N[X]$ spans all colours of~$\sigma$.

    Let us define $X_1\coloneqq \bigcup_{x\in X} N(x)$. As $G$ is locally $s$-colourable, $G[X_1]$ is decomposable into $|X|$ $(s-1)$-colourable subgraphs, so its chromatic number is at most $(s-1)|X|$. By construction, $X \setminus X_1$ is the set of isolated vertices in $G[X]$, so $G[X\setminus X_1]$ is an independent set. We infer that $\chi(G[N[X]])\le 1+(s-1)|X| \le 1+(s-1)(k+1)$.
    Since $N[X]$ spans all colours of $\sigma$, we get that 
    \[
    f_G(\chi(G) + k) \le \chi(G[N[X]]) \le (s-1)(k+1) + 1,
    \]
    as desired. 
\end{proof}

\section{Chromatic discrepancy of locally \texorpdfstring{$s$}{s}-colourable graphs}
\label{sec:discr_locally_s_colourable}

We now discuss the general case of \Cref{conjecture:main_conj}. We provide some support in favour of its statement in both directions.

\subsection{A family of locally \texorpdfstring{$s$}{s}-colourable graphs \texorpdfstring{$G$}{G} with chromatic discrepancy at most \texorpdfstring{$\chi(G)-s+1$}{χ(G)-s+1}}

In this subsection we briefly justify that, if true, \Cref{conjecture:main_conj} is almost best possible. Our construction is again based on Mycielski graphs. For integers $k\geq s \geq 2$, we recursively define $\Myc{k}{s}$ as follows. If $k=s$ then we let $\Myc{k}{s}$ be the complete graph on $s$ vertices. Otherwise, $k>s$ and $\Myc{k}{s}$ is obtained from $\Myc{k-1}{s}$ by adding a twin of each vertex and a vertex adjacent to all the twins --- that is, $\Myc{k}{s}$ is the Mycielskian of $\Myc{k-1}{s}$.

\begin{proposition}
    For all integers $k \ge s \ge 2$, $\Myc{k}{s}$ is locally $s$-colourable, has chromatic number $k$, and chromatic discrepancy $\discr(\Myc{k}{s}) \le k-s+1$. 
\end{proposition}
\begin{proof}
    It is straightforward to check by induction that $\Myc{k}{s}$ is locally $s$-colourable. It is further well-known that $\Myc{k}{s}$ has chromatic number $k$~\cite{mycielski1955}. To obtain that $\discr(\Myc{k}{s}) \leq k-s+1$, we briefly show by induction that $\Myc{k}{s}$ admits a proper $k$-colouring $\sigma_{k,s}$ such that, for any set $X\subseteq V(\Myc{k}{s})$, if $\sigma_{k,s}(X) = [s-1]$ then $X$ is a clique. It then follows that any rainbow cover $Y$ of $\sigma_{k,s}$ satisfies 
    \[ 
    |\sigma_{k,s}(Y)| - \chi(\Myc{k}{s}[Y]) \leq k-s+1,
    \]
    hence implying $\discr(\Myc{k}{s}) \leq k-s+1$, as desired.

    We define $\sigma_{k,s}$ by induction on $k$. When $k=s$ we let $\sigma_{k,s}$ be any proper $s$-colouring of $\Myc{k}{s}=K_s$ assigning to each vertex an arbitrary colour from $[s]$. It trivially satisfies the desired property.
    Assume now that $k > s$. We extend $\sigma_{k-1,s}$ to $\Myc{k}{s}$ by assigning all twins colour $k$, which is not used by $\sigma_{k-1,s}$, and colour $s$ to the vertex adjacent to all twins. By construction, every set $X$ with $\sigma_{k,s}(X) = [s-1]$ is included in $V(\Myc{k}{s})$, and by induction it follows that it is a clique.
\end{proof}

\subsection{A lower bound on the chromatic discrepancy when \texorpdfstring{$\chi$ is close to $s$}{χ is close to s}}

In this section, we show that Conjecture~\ref{conjecture:main_conj} holds when $s$ is sufficiently close to $\chi$, hence proving Theorem~\ref{thm:conj-special-case}. We also obtain a result on co-locally $(s-1)$-colourable graphs (defined later on) that we use in Section~\ref{sec:cl-free} to derive Theorem~\ref{thm:weak_conj_special_case}.

Our proof relies on a classical result on the structure of critical graphs due to Gallai.
We first introduce a few specific definition.

Given an integer $\chi\geq 1$, a graph $G$ is {\it $\chi$-critical} if $\chi(G) = \chi$ and every proper induced subgraph $H$ of $G$ satisfies $\chi(H) \leq \chi-1$.
Given two graphs $G_1$ and $G_2$, the {\it Dirac join} of $G_1$ and $G_2$, denoted $G_1 \dijoin G_2$, is the graph obtained from disjoint copies of $G_1$ and $G_2$ by adding all edges between them. 
It is straightforward that $G_1\dijoin G_2$ is $(\chi(G_1)+\chi(G_2))$-critical if and only if $G_i$ is $\chi(G_i)$-critical for $i\in \{1,2\}$.

A graph is {\it decomposable} if it is the Dirac join of two non-empty graphs, and it is {\it indecomposable} otherwise.
The following result was obtained by Gallai in 1963~\cite{gallai63b} (see~\cite{stehlikJCTB89} for a proof in english).

\begin{theorem}[Gallai~\cite{gallai63b}]
\label{thm:gallai}
    If $G$ is a $\chi$-critical indecomposable graph, then $|V(G)| \geq 2\chi-1$.
\end{theorem}

The following is a well-known consequence of Theorem~\ref{thm:gallai} (see e.g.~\cite[Corollary 4.19]{stiebitz2024}).

\begin{corollary}
    \label{cor:gallai}
    If $G$ is a $\chi$-critical graph with no universal vertex, then $|V(G)| \geq \frac{5}{3}\chi$.
\end{corollary}

From this we infer the following.

\begin{lemma}
\label{lem:size-critical}
    Let $G$ be a locally $s$-colourable graph on $n$ vertices, then 
    $\chi(G) \le \max \left\{s,  \floor{\frac{3n}{5}} \right\}$.
\end{lemma}

\begin{proof}
    Let $\chi\coloneqq \chi(G)$ and assume that $\chi>s$.
    We extract from $G$ a $\chi$-critical subgraph $H$. Since $\chi(H[N[v]])\le s < \chi$ for all vertices $v\in V(H)$, we infer that $H$ contains no universal vertex. By Corollary~\ref{cor:gallai}, $n\ge |V(H)|\ge \frac{5}{3}\chi$, so $\chi\le \frac{3}{5}n$. The chromatic number of $G$ being an integer, the result follows.
\end{proof}

A more precise application of Theorem~\ref{thm:gallai} yields the following.

\begin{lemma}
\label{lem:size-critical2}
    Let $G$ be a locally $s$-colourable graph on $n$ vertices, then 
    $\chi(G) \le \max \left\{s+1,  \floor{\frac{6n}{11}} \right\}$.
\end{lemma}

\begin{proof}
    Let us first prove by induction on $n$ that, for every locally $(\chi-2)$-colourable $\chi$-critical graph $H$, 
    \[
    |V(H)| \ge \frac{11}{6}\chi.
    \]

    Assume first that $H$ is indecomposable, hence by \Cref{thm:gallai} we have$|V(H)|\ge 2\chi-1$. Since $H$ is $\chi$-critical and locally $(\chi-2)$-colourable, we infer that $\chi\ge 4$.
    If $\chi = 4$, then $H$ is locally $2$-colourable, that is triangle-free, and so $|V(H)|\ge 11$ (see \cite{Chv06}).
    If $\chi=5$, then $H$ is locally $3$-colourable, so in particular $K_4$-free, hence $|V(H)|\ge 11$ (see \cite{Nen84}).
    The result follows since $\ceil{\frac{11k}{6}}\le 2k-1$ for every integer $k\ge 6$. 
    
    Henceforth, assume that $H = H_0 \wedge H_1$, where $H_i$ is a $\chi_i$-critical graph for each $i\in \{0,1\}$. 
    Observe that, if there exists a vertex $v\in V(H_i)$ with $\chi(H_i[N[v]])\ge \chi_i-1$, then 
    \[
    \chi(H[N_H[v]]) \ge \chi_{1-i}+\chi_i-1 = \chi-1,
    \]
    a contradiction to the fact that $H$ is locally $(\chi-2)$-colourable. So each $H_i$ is locally $(\chi_i-2)$-colourable, and by the induction hypothesis we have
    $|V(H)|\ge |V(H_1)| + |V(H_2)| \ge \frac{11}{6}(\chi_1+\chi_2) = \frac{11}{6}\chi$,
    which ends the proof of the induction.

    To finish the proof, we let $\chi\coloneqq \chi(G)$, and we assume that $\chi>s+1$.
    We extract from $G$ a locally $s$-colourable $\chi$-critical subgraph $H$ (with $s\le \chi-2$). We have $n\ge |V(H)|\ge \frac{11}{6}\chi$, so $\chi\le \frac{6n}{11}$. Since $\chi$ is an integer, the result follows.
\end{proof}

We are now able to prove \Cref{thm:conj-special-case}, that we first repeat here for convenience. 

\begin{theorem}
    For every integer $s\ge 2$, every locally $s$-colourable graph $G$ of chromatic number $\chi < \frac{11s+16}{6}$ satisfies
    \[\discr(G)\ge \chi-s.\]
\end{theorem}

\begin{proof}
    Let $p\ge \chi$ and let $\sigma \in \C_p(G)$.
    If $p = \chi$, by \Cref{lemma:rainbow_closed_neighbourhood}, there is a vertex $v \in V(G)$ such that $|\sigma(N[v])|=p=\chi$.
    Since $G$ is locally $s$-colourable, we have $\chi(G[N[v]])\le s$, and so $\discr_\sigma(G) \ge \chi-s$.    

    Henceforth, we assume that $p \ge \chi+1$.
    Let $H$ be any rainbow cover of $\sigma$, that is $|\sigma(V(H))|=|V(H)|=p$.
    By \Cref{lem:size-critical2}, we have 
    \[ \chi(H) \le \max \left\{s+1,  \floor{\frac{6p}{11}} \right\}, \]
    and so 
    \begin{align*}
        \discr_\sigma(G) &\ge p - \chi(H) \ge \min \left\{ p-s-1, \ceil{\frac{5p}{11}} \right\} \\
        &\ge  \min \left\{ \chi-s, \ceil{\frac{5(\chi+1)}{11}} \right\} \ge \chi-s,
    \end{align*}
    where we have used that $\frac{5(\chi+1)}{11} > \chi-s-1$ by the assumption that $\chi < \frac{11s+16}{6}$. This ends the proof.
\end{proof}

At the cost of more restrictive hypotheses, we can increment the lower bound in \Cref{thm:conj-special-case}. We say that a graph $G$ is \emph{co-locally $s$-colourable} for some integer $s\ge 1$ if, for every $v\in V(G)$ and every $u\in V(G) \setminus \{v\}$, $\chi(G[\{u\} \cup N(v)]) \le s$. 
Observe that every locally $s$-colourable graph is co-locally $s$-colourable. Conversely, every co-locally $s$-colourable graph is locally $(s+1)$-colourable.
We make use of the following result to prove Theorem~\ref{thm:weak_conj_special_case} in Section~\ref{sec:cl-free}.

\begin{theorem}
    \label{thm:strong-conj-special-case}
    For every integer $s\ge 3$, every co-locally $(s-1)$-colourable graph $G \neq K_s$ of chromatic number $\chi < \frac{5s+2}{3}$ satisfies
    \[ \discr(G) \ge \chi - s + 1. \]
\end{theorem}

\begin{proof}
    Since $G\neq K_s$ and $G$ is locally $s$-colourable, we infer that $G$ is not a complete graph.
    
    Let $p\ge \chi$, and let $\sigma \in \C_p(G)$.
    If $p = \chi$, since $G$ is not a complete graph, there is a colour class $U$ of size at least $2$ in $\sigma$. By \Cref{lemma:rainbow_closed_neighbourhood}, there is a vertex $v \in U$ such that $|\sigma(N[v])|=p=\chi$.
    Let $u \in U \setminus \{v\}$, and let $Y \coloneqq \{u\} \cup N(v)$. By choice of $u$, we have $|\sigma(Y)|=\chi$. Since $G$ is co-locally $(s-1)$-colourable, $\chi(G[Y]) \le s-1$. We conclude that $\discr_\sigma(G) \ge \chi-s+1$. 

    We now assume that $p \ge \chi+1$.
    Let $H$ be a rainbow cover of $\sigma$, that is $|\sigma(V(H))|=|V(H)|=p$.
    Since $H$ is locally $s$-colourable, by \Cref{lem:size-critical} we have 
    \[ \chi(H) \le \max \left\{s,  \floor{\frac{3p}{5}} \right\}. \]
    Therefore,
    \begin{align*}
        \discr_\sigma(G) &\ge p - \chi(H) \ge \min \left\{ p-s, \ceil{\frac{2p}{5}} \right\} \\
        &\ge  \min \left\{ \chi-s+1, \ceil{\frac{2(\chi+1)}{5}} \right\} \ge \chi-s+1,
    \end{align*}
    where we have used that $\frac{2(\chi+1)}{5} > \chi-s$ by the assumption that $\chi < \frac{5s+2}{3}$.
    This ends the proof.
\end{proof}

\subsection{A general lower bound on the chromatic discrepancy of locally \texorpdfstring{$s$}{s}-colourable graphs}
We now provide a short proof of Theorem~\ref{thm:locally-s-chromatic}, that we first recall here for convenience.

\begin{theorem}
    For every integer $s\ge 3$, every locally $s$-colourable graph $G$ satisfies
    \[\discr(G)\ge \chi(G)-s \ln \chi(G).\]
\end{theorem}

\begin{proof}
    Let $G$ be a locally $s$-colourable graph with chromatic number $\chi$. Let $g\colon x \mapsto x - s\ln x$.
    We prove that $\discr_\sigma(G) \geq g(\chi)$ for any proper colouring $\sigma$ of $G$. Let us thus fix such a colouring $\sigma$. 
    We proceed by induction on $\chi$. 
    
    We first briefly justify that $G$ has a $\sigma$-rainbow independent set $I$ of size at least $\chi/s$.
    We construct $I$ inductively by choosing an arbitrary vertex $v\in V(G)$, removing $N(v)$ from $G$ as well as the colour class of $v$ to obtain a graph $G'$ and adding $v$ to the rainbow independent set of $G'$ of size at least $\frac{\chi(G')}{s} \ge \frac{\chi}{s}-1$ obtained inductively. 

    If $\chi \le 3s/2$, then $g(\chi)<0$ and the result is trivial. Let us thus assume that $\chi \ge 3s/2$, and let $I$ be a $\sigma$-rainbow independent set of $G$ of size exactly $\ceil{\frac{\chi}{s}}$. Let $G'$ be obtained by removing all colour classes of $\sigma$ that intersect $I$. Observe that any proper $k$-colouring of $G'$ can be extended into a proper $(k+|I|)$-colouring of $G$ by choosing one colour for each colour class intersecting $I$. Therefore, $\chi(G) \leq \chi(G') +|I|$, from which we derive
    \begin{equation}
        \label{eq:chi-lowerbound}
        \chi(G') \geq \chi(G) - |I| \geq \chi - \ceil{\frac{\chi}{s}} = \floor{\chi(1-1/s)} \ge s, 
    \end{equation}
    where in the last inequality we used that $\chi\geq 3s/2$ and that $s\geq 3$. In particular, since $g$ is non-decreasing over $[s,+\infty)$, it follows that $g(G') \geq g(\chi(G) - |I|)$.
    
    Let $\sigma'$ be the restriction of $\sigma$ to $G'$. 
    Any graph $H \induced G'$ that realises $\discr_{\sigma'}(G')$ uses colours distinct from $\sigma(I)$, and adding $I$ to $H$ increases its chromatic number by at most $1$. We infer that
    \begin{align*}
        \discr_\sigma(G) &\ge |I| - 1 + \discr_{\sigma'}(G')\\
        &\ge |I| - 1 + g(\chi(G')) & \mbox{by the induction hypothesis;}\\
        &\ge |I| - 1 + g\pth{\chi - |I|} & \mbox{by \eqref{eq:chi-lowerbound}, since $g$ is non-decreasing over $[s,+\infty)$;}\\
        & \ge \chi - s \ln \pth{\pth{1-\frac{1}{s}}\chi} - 1 \\
        & = \chi -s \ln \pth{\chi} - s \ln \pth{1-\frac{1}{s}}  - 1 \\
        & \ge g(\chi) &\mbox{since $-\ln(1-x) \ge x$ for all $0\le x < 1$.}
    \end{align*}
    This ends the proof of the induction.
\end{proof}

\section{Chromatic discrepancy of \texorpdfstring{$2$}{2}-locally \texorpdfstring{$s$}{s}-colourable graphs}
\label{sec:2-local}

The goal of this section is to prove Theorem~\ref{thm:2-local}, which states that every $1$-locally $s_1$-colourable and $2$-locally $s_2$-colourable graph $G$ satisfies $\discr(G) \geq \chi(G) - O(s_1s_2\ln \ln(\chi(G)))$. 
We first overview the proof briefly. Assume that $G$ is such a graph with chromatic number $\chi$, and consider any proper colouring $\sigma$ of $G$. We divide the proof into two cases, depending on whether the number of colours $p$ used by $\sigma$ is close enough to $\chi$ or not. 

The easiest case is when $p$ is large enough compared to $\chi$, that is $p \geq \chi + \Omega(\sqrt{\chi})$. In that case, we prove that any $\sigma$-rainbow cover $X$ satisfies $|\sigma(X)| - \chi(G[X]) \geq \chi$, hence implying the result.
When $p$ is close to $\chi$, that is $p \leq \chi + O(\sqrt{\chi})$, with Lemma~\ref{lem:extract} in hand we manage to find an independent set $I$ that covers almost all colours of $\sigma$. Formally, $I$ is an independent set such that $|\sigma(X)| \geq p - p^{1-\varepsilon}$, where the value of $\varepsilon>0$ depends only on $s_1,s_2$. Removing all colour classes intersecting $I$ and adding $I$ to a rainbow set $X'$ with small chromatic number obtained inductively yields a $\sigma$-rainbow cover $X$ with small chromatic number again, hence showing $\sigma(X) - \chi(G[X]) \geq \chi - O(s_1s_2\ln \ln(\chi))$.

\medskip

We start with the case of $p$ being large enough compared to $\chi$ that, in the proof of Theorem~\ref{thm:2-local}, we cover with the following lemma.

\begin{lemma}
    \label{lem:chi_1local}
    Let $G$ be a $1$-locally $s$-colourable graph on $n$ vertices, then $\chi(G) \le \sqrt{2sn}$.
\end{lemma}

\begin{proof}
    We prove the result by induction on $n$. 
    If $n \le 2s$, then the bound is trivial. 
    Let us now assume that $n>2s$. 
    If $\Delta(G) \le \sqrt{2sn}-1$, then the result holds directly, as an easy greedy procedure shows that $\chi(G)\leq \Delta(G)+1$.
    Otherwise, there is a vertex $v$ of degree at least $\sqrt{2sn}$ in $G$. Since $G[N[v]]$ is $s$-colourable, it contains an independent set $I$ of size at least $\sqrt{2n/s}$ by the Pigeonhole Principle. 
    We apply induction on $G_0 \coloneqq G\setminus I$, and use the bound $f(x-h) \le f(x) - h f'(x)$ for every function $f$ concave on $[x-h,x]$. We obtain
    \begin{align*}
        \chi(G_0) &\le \sqrt{2s(n-|I|)} \le \sqrt{2sn} - |I| \sqrt{\frac{s}{2n}} &\mbox{by concavity of the function $x\mapsto \sqrt{2sx}$}\\
        &\le \sqrt{2sn} - 1.
    \end{align*}
    Since $G_0$ is obtained by removing an independent set from $G$, we have $\chi(G) \le \chi(G_0)+1 \le \sqrt{2sn}$, as desired.
\end{proof}

We note that one could easily improve the bound of \Cref{lem:chi_1local} to one of the form $O_s(\sqrt{n/\log n})$ by relying on stronger bounds than the greedy one for the chromatic number of locally $s$-colourable graphs of maximum degree $\Delta$ (see e.g. \cite{DKPS20+}). However, such a level of precision is not needed in our proof, so we prefer this easier bound with a self-contained proof.

We now consider the other case, that is $p$ being close to $\chi$. Our goal here is to show the existence of an independent set spanning almost all $p$ colours. To show the existence of such an independent set, we make use of the following direct consequence of Lemma~\ref{lemma:rainbow_closed_neighbourhood} via an application of the Pigeonhole Principle.

\begin{lemma}\label{lem:extract}
    Let $G$ be a graph, and let $p=\chi(G)+k$ for some integer $k$. Then, for every $\sigma\in \C_p(G)$, there is a vertex $v\in V(G)$ whose closed neighbourhood contains at least $\frac{p}{k+1}$
    colours in $\sigma$.
\end{lemma}

\begin{lemma}
\label{lem:RIS-2-local}
    For all integers $s_1,s_2\ge2$, the following holds.
    Let $G$ be a $1$-locally $s_1$-colourable and $2$-locally $s_2$-colourable graph of chromatic number $\chi$, let $p=\chi+k$ be an integer with $k\ge 0$, and let $\sigma \in \C_p(G)$. Then $G$ has a $\sigma$-rainbow independent set of size at least 
    \[p-p^{1-\frac{1}{3s_1s_2}} - \frac{k+1}{s_2}p^{1/3}.\]
\end{lemma}

\begin{proof}
    Before delving into the formal proof, let us describe its key idea. If a graph $G$ is properly coloured using not much more than $\chi(G)$ colours, \Cref{lemma:rainbow_closed_neighbourhood} ensures that there is a vertex $v\in V(G)$ whose neighbourhood spans a non-negligible fraction of those colours, and if moreover $N(v)$ has bounded chromatic number --- this is ensured by the fact that $G$ is $1$-locally $s_1$-colourable ---, we can extract from $N[v]$ a large rainbow independent set $I$. One can further extend $I$ by repeating the argument on the graph obtained from $G$ by removing $N(I)$ and all colour classes that intersect $I$. By doing so, the difference $k$ between the total number of colours used and the chromatic number does not increase too much, as long as $N(I)$ has bounded chromatic number --- this is ensured by the fact that $G$ is $2$-locally $s_2$-colourable since $N(I) \subseteq N^2[v]$. We now describe that process formally, by carefully keeping track of all parameters at play.

    Let $G$ be a $1$-locally $s_1$-colourable and $2$-locally $s_2$-colourable graph of chromatic number $\chi$ and let $\sigma\in \C_p(G)$ be a proper $p$-colouring of $G$ with $p=\chi+k$.
    Let $G_0\coloneqq G$, $p_0\coloneqq p$, and $I_0 \coloneqq \emptyset$.
    We recursively define $(I_i,G_i,p_i)$ for all integers $i\in \mathbb{N}$, maintaining the following properties:
    \begin{enumerate}[label=(\roman*)]
        \item $I_i$ is a rainbow independent set of $G$ of size $p-p_i$;
        \item $G_i$ is an induced subgraph of $G$ such that $\sigma(V(G_i)) \cap \sigma(I_i) = \emptyset$ and $V(G_i)\cap N_G(I_i) = \emptyset$ --- in particular, $|\sigma(V(G_i))|\le p_i$;
        \item $k_i \coloneqq p_i - \chi(G_i) \le k_0 + is_2$.
    \end{enumerate}
    While $G_i$ is non-empty, we define $(p_{i+1},I_{i+1},G_{i+1})$ from $(p_i,I_i,G_i)$ as follows.
    \begin{enumerate}[label=(\arabic*)]
        \item We set $p_{i+1} \coloneqq p_i - \floor{\frac{p_i}{s_1(k_i+1)}}$.
        \item Let $v\in V(G_i)$ be such that 
        \[
        |\sigma(N_{G_i}[v])|\ge \frac{|\sigma(V(G_i))|}{|\sigma(V(G_i))|-\chi(G_i)+1},
        \]
        the existence of which is guaranteed by Lemma~\ref{lem:extract}.
        Observe that, when $G_i$ is non-empty and thus $\chi(G_i)\geq 1$, the function $h\colon x\mapsto \frac{x}{x-\chi(G_i)+1}$ is non-increasing on $[1,+\infty)$. Therefore, since $|\sigma(V(G_i))|\leq p_i$, we further have
        \[
        |\sigma(N_{G_i}[v])|\ge \frac{|\sigma(V(G_i))|}{|\sigma(V(G_i))|-\chi(G_i)+1} \ge \frac{p_i}{p_i-\chi(G_i)+1} = \frac{p_i}{k_i+1}.
        \]
        As $G_i$ is a subgraph of $G$, we have that $G_i[N[v]]$ is $s_1$-colourable, meaning that we can extract in $G_i[N[v]]$ a rainbow independent set $J_i$ of size exactly 
        \[ |J_i| = \floor{\frac{p_i}{s_1(k_i+1)}}.\]
        It follows from (ii) that $I_{i+1}\coloneqq I_i \cup J_i$ is a rainbow independent set of $G$ of size at least $p-p_i + |J_i| = p - p_{i+1}$ and (i) is maintained.
        \item Let $X_i$ be the union of the colour classes of $\sigma$ that intersect $J_i$. 
        We set $G_{i+1} \coloneqq G_i \setminus (X_i \cup N(J_i))$, which maintains (ii).

        Observe that $N(J_i)$ is contained in the second neighbourhood of $v$, so, since $G$ is $2$-locally $s_2$-colourable, $\chi(G[N(J_i)])\le s_2$. Moreover, $\chi(G[X_i])\leq |J_i|$ as $\sigma$ yields a proper $|J_i|$-colouring of $G[X_i]$ using exactly $|J_i|$ colours. 
        Therefore, $\chi(G_{i+1}) \ge \chi(G_i) - |J_i| - s_2$, hence
        \begin{align*}
            k_{i+1} = p_{i+1} - \chi(G_{i+1}) &\le p_{i+1} - \chi(G_i) + |J_i| + s_2 \\ 
            &= k_i + s_2 \\
            &\le k_0 + (i+1)s_2,
        \end{align*}
        and (iii) is maintained.
    \end{enumerate}
If $G_i$ is empty, we let $(p_{i+1},I_{i+1},G_{i+1})\coloneqq (p_i,I_i,G_i)$. Note that, for such $i$, one has $k_i=p_i$, and so $\floor{\frac{p_i}{s_1(k_i+1)}}=0$. Therefore, for all integers $i\ge 0$, one has
\begin{equation}
    \label{eq:p_i}
    p_{i+1} = p_i - \floor{\frac{p_i}{s_1(k_i+1)}} \le p_i\pth{1 - \frac{1}{s_1(k_i+1)}}+1 \le p_i e^{- \frac{1}{s_1(k_i+1)}}+1.
\end{equation}
We now provide an estimation of the value of $p_i$. Using the well-known fact that $\sum_{j=0}^{i-1} g(x) \ge \int_0^i g(x)dx $ for every non-increasing function $g$ over the interval $[0,i]$, we derive from \eqref{eq:p_i} that
\begin{align*}
    p_i &\le  p \exp\pth{- \sum_{j=0}^{i-1} \frac{1}{s_1(k_j+1)}} + i \le p \exp\pth{- \sum_{j=0}^{i-1} \frac{1}{s_1(k_0 + js_2+1)}} + i \\
    &\le p \exp\pth{- \int_{0}^{i} \limits \frac{dx}{s_1(k_0 + s_2x+1)}} + i = p \exp\pth{-\frac{1}{s_1s_2} \ln \pth{\frac{k_0 + is_2+1}{k_0+1}}} + i.
\end{align*}
By fixing $i \coloneqq \frac{k_0+1}{s_2}p^{1/3}$, we obtain that $I_i$ is a rainbow independent set of size
\[ 
p-p_i \ge p - p^{1-\frac{1}{3s_1s_2}} - \frac{k_0+1}{s_2}p^{1/3},
\]
as desired.
\end{proof}

We are now ready to prove the main result of this section, namely Theorem~\ref{thm:2-local}, whose proof relies on a combination of Lemmas~\ref{lem:chi_1local} and~\ref{lem:RIS-2-local}. We first recall it for convenience.

\begin{theorem}
    For all integers $s_1,s_2 \ge 2$, there exists an absolute constant $C\ge 0$ such that the following holds.
    Every $1$-locally $s_1$-colourable and $2$-locally $s_2$-colourable graph $G$ satisfies
    \[
    \discr(G) \ge \chi - 4s_1s_2\ln \ln \chi - C.
    \]
\end{theorem}

\begin{proof}
    We do not provide the explicit value of $C$ and simply assume that it is large enough in terms of $s_1$ and $s_2$. 
    Let $g\colon x \mapsto x-4s_1s_2\ln \ln x - C$ and observe that this function is non-decreasing over the interval $[4s_1s_2, +\infty)$. Indeed, the derivative of $g$ is $g'(x) = 1 - \frac{4s_1s_2}{x\ln x} \ge 0$ when $x \ge 4s_1s_2$.
    
    We now prove by induction on $|V(G)|$ that $\discr(G) \ge g(\chi(G))$. Note that the result is trivial when $\chi \le C$, since then $g(\chi)\le 0$. We thus assume that $\chi \ge C$. 
    Let $\sigma$ be any proper $p$-colouring of $G$ such that
    \[
        \discr(G) = \discr_\sigma(G) = p-f_G(p),
    \]
    the existence of which is guaranteed by Lemma~\ref{lem:discrfp}.
    We now distinguish two cases, depending on the value of $p$ with respect to $\chi(G)=\chi$. 
    
    Assume first that $p \geq \chi + 2s_2\sqrt{\chi}$. Let $X$ be any rainbow cover of $\sigma$. Observe that any $2$-locally $s$-colourable graph is $1$-locally $s$-colourable graph by definition. By \Cref{lem:chi_1local}, $G[X]$ has thus chromatic number at most $\sqrt{2s_2p}$. Hence, $\discr(G) \geq p-\sqrt{2s_2p}$.
    Note that, $C$ being large enough, the function $x\mapsto x-\sqrt{2s_2x}$ is non-decreasing on $[C,+\infty)$. Therefore, it follows from the inequality above and the fact that $p \geq \chi + 2s_2\sqrt{\chi} \geq C$ that
    \[
        \discr(G) \geq (\chi+2s_2\sqrt{\chi})-\sqrt{2s_2(\chi+2s_2\sqrt{\chi})} > \chi,
    \]
    where, in the last inequality, we used that $\chi$ is large enough. This implies the result.

    Henceforth, assume that $p \leq \chi + 2s_2\sqrt{\chi}$.
    By \Cref{lem:RIS-2-local}, there is a $\sigma$-rainbow independent set $I$ of $G$ of size at least 
    \begin{align*}
        p-p^{1-\frac{1}{3s_1s_2}}-\frac{p-\chi+1}{s_2}p^{1/3} &\geq p-p^{1-\frac{1}{3s_1s_2}}-\frac{2s_2\sqrt{\chi}+1}{s_2}p^{1/3}\\ &\geq p-p^{1-\frac{1}{3s_1s_2}}-2p^{5/6} - p^{1/3}\\
        &\geq p - 3 p^{1-\frac{1}{3s_1s_2}},
    \end{align*}
    where in the last inequality we use that $p$ is large enough. We let $I$ be a $\sigma$-rainbow independent set of size exactly $\ceil{p - 3 p^{1-\frac{1}{3s_1s_2}}}$.
        
    Let $G'$ be obtained from $G$ by removing all the colour classes of $\sigma$ that intersect $I$. Assuming that $C$, and therefore also $\chi$ and $p$, is large enough, we have 
    \begin{equation}
    \label{eq:chi'-lower}
        \chi(G') \ge \chi-|I| = \floor{\chi - (p-3p^{1-\frac{1}{3s_1s_2}})} \geq \floor{3(\chi + 2s_2\sqrt\chi)^{1-\frac{1}{3s_1s_2}} - 2s_2\sqrt\chi} \ge 4s_1s_2,
    \end{equation}
    where in the third inequality we use that $h\colon x\mapsto x-3x^{1-\frac{1}{3s_1s_2}}$ is non-decreasing over $[C,+\infty)$ and that $p\leq \chi+2s_2\sqrt\chi$. In particular,~\eqref{eq:chi'-lower} implies that $g(\chi(G')) \geq g(\chi - |I|)$, since $g$ is non-decreasing over $[4s_1s_2,+\infty)$. Using again that $\chi$ and $p$ are large, we also obtain
    \begin{equation}
        \label{eq:p-lower}
        3p^{1-\frac{1}{3s_1s_2}} \leq 3(\chi +2s_2\sqrt\chi)^{1-\frac{1}{3s_1s_2}}\leq  \chi^{1-\frac{1}{4s_1s_2}}.
    \end{equation}
    We infer that
    \begin{align*}
        \discr_\sigma(G) &\ge |I|-1 + \discr(G') \\
        &\ge |I|-1 + g(\chi(G')) & \mbox{by the induction hypothesis;}\\
        &\ge |I|-1 + g(\chi-|I|) & \mbox{by \eqref{eq:chi'-lower};}\\
        &= \chi - 1 - 4s_1s_2 \ln \ln \pth{\chi - |I|} - C\\
        &\ge \chi - 1 - 4s_1s_2 \ln \ln \pth{3p^{1-\frac{1}{3s_1s_2}}} - C & \mbox{since $p\geq \chi$;}\\
        &\ge \chi - 1 - 4s_1s_2 \ln \ln \pth{\chi^{1-\frac{1}{4s_1s_2}}} - C & \mbox{by \eqref{eq:p-lower}}\\
        &= \chi - 1 - 4s_1s_2 \pth{\ln\pth{1-\frac{1}{4s_1s_2}} + \ln \ln \chi} - C \\
        &\ge \chi - 4s_1s_2 \ln \ln \chi -C &\mbox{since $\ln(1+x) \le x$ for all $x\ge -1$.}
    \end{align*}
    The result follows.
\end{proof}

\section{Chromatic discrepancy of graphs excluding a cycle}
\label{sec:cl-free}

Recall that $C_{\ell+1}$-free graphs are locally $\ell$-colourable by Proposition~\ref{prop:C_ell_Nv_degenerate}. The following is thus a direct consequence of Theorem~\ref{thm:chi_plus_k}, where $f_G$ is the function defined in \eqref{eq:def-f_G}.

\begin{corollary}
    \label{cor:chi+k}
    Let $\ell\ge 2$, $k\ge 0$ be two integers. Every $C_{\ell+1}$-free graph $G$ satisfies \[f_G(\chi(G) + k) \le (\ell-1)(k+1)+1.\]
\end{corollary}

Unfortunately, the bound given by Corollary~\ref{cor:chi+k} is not strong enough to deduce, using Lemma~\ref{lem:discrfp}, a lower bound on $\discr(G)$. However, we may use the finer structure exhibited by Lemma~\ref{lemma:rainbow_closed_neighbourhood} to improve Corollary~\ref{cor:chi+k} for $C_4$-free graphs, hence showing Theorem~\ref{thm:C4_free}. We first recall it for convenience. 

\begin{theorem}
    Every $C_4$-free graph $G\neq K_3$ satisfies
    \[\discr(G)\ge \chi(G)-2.\]
\end{theorem}

\begin{proof}
    Let $\chi=\chi(G)$. We assume that $\chi\geq 3$, the result being trivial otherwise.
    We claim that $f_G(\chi+k) \le k+2$ for every $k\ge 0$. By Lemma~\ref{lem:discrfp}, this implies 
    \[ \discr(G)\ge \min \sset{\chi+k-(k+2)}{k\ge 0}=\chi-2,\]
    as desired.
    To prove the claim, we fix an integer $k\ge 0$, and we let $\sigma \in \C_{\chi+k}(G)$ be a proper colouring of $G$ using $\chi+k$ colours. 

    We first treat the case $k=0$ separately.
    Since $G\neq K_3$ and $G$ is $C_4$-free, we infer that $G$ is not a complete graph.
    In particular, $\chi(G) < |V(G)|$, and there is a colour class $U$ of size at least $2$ in $\sigma$. 
    By \Cref{lemma:rainbow_closed_neighbourhood}, there is a vertex $v \in U$ such that $|\sigma(N[v])|=p=\chi$.
    Since $G$ is $C_4$-free, by Proposition~\ref{prop:C_ell_Nv_degenerate} we have $\chi(G[N(v)])\le 2$. 
    Let $u \in U \setminus \{v\}$, and let $Y \coloneqq \{u\} \cup N(v)$, so that $|\sigma(Y)|=\chi$. Since $G$ is $C_4$-free, $u$ has at most one neighbour in $N(v)$, so $\chi(G[Y]) \le 2$. We conclude that $f_G(\chi) \le 2$. 

    We now assume that $k\ge 1$. 
    By Lemma~\ref{lemma:rainbow_closed_neighbourhood}, there exists a set $X$ of size at most $k+1$ such that $\sigma(N[X]) = [\chi + k]$. 
    We show that $G[N[X]]$ is $(k+2)$-colourable, which implies $f_G(\chi+k)\le k+2$, as desired.
    
    Let $H \coloneqq G[N(X) \setminus X]$.
    For every $u\in V(H)$, let $h(u)$ be a fixed neighbour of $u$ in $X$. We let $Y_i \coloneqq \sset{u\in V(H)}{h(u)=x_i}$ for each $i\in [k+1]$. Observe that if a vertex $u\in V(H)$ has two neighbours $u_1,u_2\in Y_i$ for some $i$, then $u,u_1,x_i,u_2$ is a $4$-cycle, a contradiction. So every vertex $u\in V(H)$ has at most one neighbour in each $Y_i$. 
    
    We first treat the case of $X$ being an independent set of $G$. In this case, we show that $\chi(H)\le k+1$, which implies that $\chi(G[N[X]])\le k+2$ since $X$ is an independent set. 
    If $k=1$, we colour each edge in $H$ blue if its endpoints lie in the same set $Y_i$, or red otherwise. This is a proper $2$-edge-colouring of $H$, witnessing that $H$ has no odd cycle, and that $\chi(H)\le 2$, as desired.
    If $k\ge 2$, we observe that $\Delta(H) \le k+1$, and since $\omega(H)\le 3 < k+2$, we infer by Brooks' Theorem (Theorem~\ref{thm:brooks}) that $\chi(H) \le k+1$, as desired.

    Henceforth, we assume that $X$ is not an independent set. 
    Let $Z \subseteq N(X)\setminus X$ be the set of vertices outside of $X$ with at least $2$ neighbours in $X$. We first prove that $G[X\cup Z]$ is $(k+2)$-colourable.
    \begin{claim}
        $\chi(G[X\cup Z]) \leq k+2$.
    \end{claim}
    \begin{proofclaim}
        Assume for the sake of contradiction that $\chi(G[X\cup Z])\ge k+3$, so in particular there is a subgraph $H_0 \subseteq G[X\cup Z]$ of minimum degree at least $k+2$. It follows that $|V(H_0)| \ge k+3$ and that $|V(H_0)\cap Z| \ge 2$. 
        
        For every $z\in Z$, observe that, if there exist $z_1 \neq z_2 \in N_Z(z)$, then $N_X(z_1)\cap N_X(z_2) = \emptyset$, otherwise there would be a $4$-cycle that contains $z_1,z,z_2$ and a common neighbour of $z_1$ and $z_2$ in $X$.
        Since every vertex in $Z$ has two neighbours in $X$ by definition, by the Pigeonhole Principle, this implies $\deg_Z(z) \le |X|/2 \le \frac{k+1}{2}$.
        We infer that, for every vertex $z \in V(H_0) \cap Z$, 
        \[\deg_X(z) = \deg_{H_0}(z) - \deg_Z(z) \ge k+2 - \frac{k+1}{2} = \frac{k+3}{2},\]
        where we have used that  $H_0$ has minimum degree at least $k+2$. Since there exist two such vertices, by the Pigeonhole Principle, they have at least $2$ common neighbours in $X$, which contradicts the fact that $G$ is $C_4$-free.
    \end{proofclaim}

    Let us thus fix a proper $(k+2)$-colouring $\sigma$ of $G[X\cup Z]$. We now show that one can extend $\sigma$ to $N(X) \setminus (X \cup Z)$. 
    Since $X$ is not an independent set, without loss of generality, let us assume that $x_1x_2\in E(G)$. 
    Since $G$ is $C_4$-free, there is no edge between $Y_1$ and $Y_2$. By definition of $Z$, it follows that $\deg_{G[N[X]]}(u) \le k+1$ for every $u\in (Y_1 \cup Y_2)\setminus Z$. 
    
    Let $K$ be a connected component of $G[N(X)\setminus (X\cup Z)]$. If $K$ does not intersect $Y_1\cup Y_2$, then for every $u\in K$, $\deg_{G[N[X]]}(u) \le k$. We conclude that each connected component $K$ of $G[N(X)\setminus (X\cup Z)]$ contains at least one vertex $v_K$ of degree at most $k+1$ in $N[X]$, and all other vertices have degree at most $k+2$ in $N[X]$. 
    
    We can therefore extend $\sigma$ to each connected component $K$ of $G[N(X)\setminus (X\cup Z)]$ as follows. Take a spanning tree $T_K$ rooted in $v_K$, and greedily colour the vertices, following a leaves-to-root ordering of $T_K$. At each step, if the considered vertex is distinct from $v_K$, then it has degree at most $k+2$ and an uncoloured neighbour (namely its father in the tree), so one of the $k+2$ colours is available. At the very last step, we end up with $v_K$, which has degree $k+1$ at most, and can therefore be coloured.
\end{proof}

Recall that a graph $G$ is co-locally $s$-colourable if $\chi(G[\{u\}\cup N(v)])\le s$ for every $\{u,v\}\in \binom{V(G)}{2}$. 

\begin{lemma}
    \label{lem:co-locally-c_ell-free}
    Let $\ell \ge 3$ be an integer.
    If a graph $G$ is $C_{\ell+1}$-free, then it is co-locally $(\ell-1)$-colourable.
\end{lemma}

\begin{proof}
    Let $G$ be a graph that is not co-locally $(\ell-1)$-colourable, i.e. there exists $\{u,v\}\in \binom{V(G)}{2}$ such that $\chi(G[Y])\ge \ell$, where $Y\coloneqq \{u\}\cup N(v)$. In particular, $G[Y]$ contains a subgraph $H$ of minimum degree at least $\ell-1\ge 2$. Let $\{x,y\} \subseteq N_H(u) \subseteq N_G(v)$; we extend greedily a path starting with $x,u,y$ into a copy of $P_\ell$ in $H$. Together with $\{v\}$, this forms a copy of $C_{\ell+1}$, so $G$ is not $C_{\ell+1}$-free, as desired.
\end{proof}

As a consequence of \Cref{lem:co-locally-c_ell-free} together with \Cref{thm:strong-conj-special-case}, we obtain the following, which corresponds to Theorem~\ref{thm:weak_conj_special_case}.

\begin{corollary}
    For every integer $\ell \ge 3$, every $C_{\ell+1}$-free graph $G\neq K_\ell$ with $\chi(G) < \frac{5\ell+2}{3}$ satisfies
    \[ 
    \discr(G) \ge \chi - \ell + 1. 
    \]
\end{corollary}

\section{Balls of small radius have bounded chromatic number in \texorpdfstring{$C_{\ell+1}$-free graphs}{graphs with a forbidden cycle length}.}
\label{sec:coloring-balls}
In this section, we show that $C_{\ell+1}$-free graphs fall into the scope of the results from Section~\ref{sec:2-local}. 
Recall that, by Proposition~\ref{prop:C_ell_Nv_degenerate}, $C_{\ell+1}$-free graphs $G$ satisfy $\chi(B_1(v))\le \ell$ for every vertex $v\in V(G)$. The following theorem is a qualitative improvement of this result, since its shows a linear bound on the chromatic number of the balls of radii $\frac{\ell}{2}$.

\begin{theorem}
\label{thm:cl}
    Let $\ell \ge 2$ be a fixed integer, and let $G$ be a $C_{\ell+1}$-free graph. We denote $t \coloneqq \floor{\frac{\ell}{2}}$. Then, for every vertex $v\in V(G)$,
    \[ \chi(B_t(v)) \le 2\ell.\]
\end{theorem}

In particular, it follows from \Cref{thm:cl} that, for every $\ell\geq 4$, $C_{\ell+1}$-free graphs are $2$-locally $2\ell$-colourable. Since they are also $1$-locally $\ell$-colourable by Proposition~\ref{prop:C_ell_Nv_degenerate}, \Cref{thm:2-local} implies the following result as a direct corollary, and further shows Theorem~\ref{thm:discr_c_ell}.

\begin{corollary}\label{cor:phi_cl}
    For every integer $\ell \ge 4$, there is a constant $c_\ell$ such that the following holds.
    For every $C_{\ell+1}$-free graph $G$ of chromatic number $\chi$,
    \[\discr(G) \ge \chi - 8\ell^2 \ln \ln \chi - c_\ell.\]
\end{corollary}

\Cref{thm:cl} follows from the following lemma, which has a different proof depending on the parity of $\ell$. 
    
\begin{lemma}
\label{lem:Lr(v)}
    Let $\ell\ge 2$ be a fixed integer, and let $G$ be a $C_{\ell+1}$-free graph. Then for every vertex $v\in V(G)$ and every integer $1\le r\le \floor{\ell/2}$, the subgraph $L_r(v)$ is $(\ell-1)$-degenerate.
\end{lemma}

Before going on with the proof of \Cref{lem:Lr(v)} in two different subsections, depending on the parity of $\ell$, we show how to derive \Cref{thm:cl} from it.

\begin{proof}[Proof of \Cref{thm:cl}]
    Let $G$ be a $C_{\ell+1}$-free graph, and let $v\in V(G)$ be any vertex. Let us write $t\coloneqq \floor{\ell/2}$.
    By Lemma~\ref{lem:Lr(v)}, for every $r\le t$, the subgraph $L_r(v)$ is $(\ell-1)$-degenerate, hence $\ell$-colourable. The statement follows, since for every $r\neq r'$ of the same parity, there is no edge between $N^r(v)$ and $N^{r'}(v)$. In particular, we may colour $\bigcup_{r \text{ even}} N^r(v)$ with $\ell$ colours, and $\bigcup_{r' \text{ odd}} N^{r'}(v)$ with $\ell$ other colours, thus yielding a proper $2\ell$-colouring of $B_t(v)$.
\end{proof}

\subsection{Forbidding an odd cycle}
We prove the following, which corresponds to the case of \Cref{lem:Lr(v)} where we forbid an odd cycle.

\begin{lemma}
\label{lem:Lr(v)_odd}
    Let $t\ge 1$ be a fixed integer, and let $G$ be a $C_{2t+1}$-free graph. Then for every vertex $v\in V(G)$ and every integer $1\le r\le t$, the subgraph $L_r(v)$ is $(2t-1)$-degenerate.
\end{lemma}
\begin{proof}
    We proceed by induction on $r$.
    When $r=1$, $L_1(v)$ is $(2t-1)$-degenerate by Proposition~\ref{prop:C_ell_Nv_degenerate}.

    Assume now that $r\ge 2$. Let $T$ be a BFS-tree of depth $r$ rooted in $v$. We let $u_1, \ldots, u_d$ be the neighbours of $v$, and $T_i$ be the subtree of $T$ rooted in $u_i$ for each $i\in [d]$. 
    Let us assume for the sake of contradiction that there is a subgraph $H\subseteq L_r(v)$ of minimum degree $\delta(H)\ge 2t$.     
    We denote $X_i \coloneqq V(H) \cap V(T_i)$ the set of descendants of $u_i$ in $V(H)$, for every $i\in [d]$. 

    \begin{claim}\label{claim:neighbourhood_Xj}
        For every $w\in V(H)$, there is $j\in [d]$ such that $N_H(w) \subseteq X_j$.
    \end{claim}
    \begin{proofclaim}
        Assume for the sake of contradiction that $w$ has a neighbour $x_i \in X_i$ and another $x_j \in X_j$, with $i\neq j$. See Figure~\ref{fig:construction_Q1_Q2_P} for an illustration.
        By the minimum degree condition on $H$, one can construct greedily a path $P$ of length $2t-2r$ starting from $w$ in $H$, and disjoint from $\{x_i,x_j\}$. \Wlog, by symmetry of the roles of $i$ and $j$, we may assume that the last vertex $z$ of $P$ belongs to $X_k$ for some $k\neq i$.

        \begin{figure}
            \centering
            \begin{tikzpicture}
            \tikzset{vertex/.style = {circle,fill=black,minimum size=5pt, inner sep=0pt}}
            \tikzset{edge/.style = {very thick}}
            \tikzset{bigvertex/.style = {shape=circle,draw}}
            \tikzset{leaves/.style ={rectangle, draw, minimum width=1.8cm, minimum height=0.8cm}}
            \definecolor{gblue}{rgb}{0.0, 0.5, 1.0}

            \node[vertex, label=above:$v$] (v) at (0,5) {};
            \node[] at (-8,3) {$N(v)$};
            \node[vertex, label=left:$u_1$] (u1) at (-6,3) {};
            \node[vertex, label=left:$u_2$] (u2) at (-4,3) {};
            \node[] (ub1) at (-2.75,3) {$\cdots$};
            \node[vertex, red, label=left:$u_i$] (ui) at (-1.5,3) {};
            \node[] (ub2) at (-0.25,3) {$\cdots$};
            \node[vertex, label=right:$u_s$] (us) at (1,3) {};
            \node[] (ub3) at (2.25,3) {$\cdots$};
            \node[vertex, label=right:$u_k$, green!70!black] (uk) at (3.5,3) {};
            \node[] (ub4) at (4.75,3) {$\cdots$};
            \node[vertex, label=right:$u_d$] (ud) at (6,3) {};
            
            \foreach \x in {1,2,s,d}{
                \draw (v) -- (u\x);
            }
            \draw[edge, red] (v) -- (ui);
            \draw[edge, green!70!black] (v) -- (uk);
            
            \node[] at (-8,0) {$N^r(v)$};
            \node[leaves, label=below:$X_1$] (X1) at (-6,0) {};
            \node[leaves, label=below:$X_2$] (X2) at (-4,0) {};
            \node[] at (-2.75,0) {$\cdots$};
            \node[leaves, label=below:$X_i$] (Xi) at (-1.5,0) {};
            \node[] at (-0.25,0) {$\cdots$};
            \node[leaves, label=below:$X_s$] (Xs) at (1,0) {};
            \node[] at (2.25,0) {$\cdots$};
            \node[leaves, label=below:$X_k$] (Xk) at (3.5,0) {};
            \node[] at (4.75,0) {$\cdots$};
            \node[leaves, label=below:$X_d$] (Xd) at (6,0) {};
            \foreach \i in {1,2,i,s,k,d}{
                \draw (u\i) -- (X\i.north east);
                \draw (u\i) -- (X\i.north west);
            }
            
            \node[] at (-6.9,1.5) {$T_1$};
            \node[] at (-4.9,1.5) {$T_2$};
            \node[] at (-2.4,1.5) {$T_i$};
            \node[] at (1.9,1.5) {$T_s$};
            \node[] at (4.4,1.5) {$T_k$};
            \node[] at (6.9,1.5) {$T_d$};

            \node[vertex, label=left:$x_i$] (xi) at (Xi) {};
            \node[vertex, label=right:$w$] (w) at (Xs) {};
            \node[vertex, label=left:$z$] (z) at (Xk) {};
            \draw[edge] (w) to[out=160, in=20] (xi);

            \node[vertex, red] (redinter1) at (-1.6, 2) {};
            \node[vertex, red] (redinter2) at (-1.3, 1) {};
            \draw[edge, red] (ui) -- (redinter1);
            \draw[edge, red] (redinter2) -- (redinter1);
            \draw[edge, red] (xi) -- (redinter2);
            \node[] at (-0.5,3.8) {\color{red} $Q_1$}; 
            
            \node[vertex, green!70!black] (greeninter1) at (3.65, 2) {};
            \node[vertex, green!70!black] (greeninter2) at (3.4, 1) {};
            \draw[edge, green!70!black] (uk) -- (greeninter1) -- (greeninter2) -- (z); 
            \node[] at (1.4,3.8) {\color{green!70!black} $Q_2$}; 

            \node[vertex, gblue] (x1) at (-6.5,0) {};
            \node[vertex, gblue] (x1b) at (-5.5,0) {};
            \node[vertex, gblue] (x2) at (X2) {};
            \node[vertex, gblue] (xd1) at (5.5,0) {};
            \node[vertex, gblue] (xd2) at (6.5,0) {};

            \draw[edge, gblue] (w) to[out=-135,in=-45] (x1b);
            \draw[edge, gblue] (x1b) to (x2);
            \draw[edge, gblue] (x2) to[out=140, in=40] (x1);
            \draw[edge, gblue] (x1) to[out=-90,in=180] (-5,-2) to (5,-2) to[out=0,in=-90] (xd2);
            \draw[edge, gblue] (xd2) to (xd1) to[out=156, in=24] (z);

            \node[] at (-0.6,-1.5) {\color{gblue} $P$};            
            \end{tikzpicture}
            \caption{An illustration of the construction of the paths $P$, $Q_1$, and $Q_2$. The path $P$, in blue, is included in $L_r(v)$. The paths $Q_1$ and $Q_2$, in red and green respectively, intersect $N^r(v)$ exactly on $\{x_i\}$ and $\{z\}$ respectively. Both $Q_1$ and $Q_2$ have length $r$. The subtrees $T_i$ being pairwise disjoint, $Q_1$ and $Q_2$ intersect exactly on $\{v\}$.}
            \label{fig:construction_Q1_Q2_P}
        \end{figure}
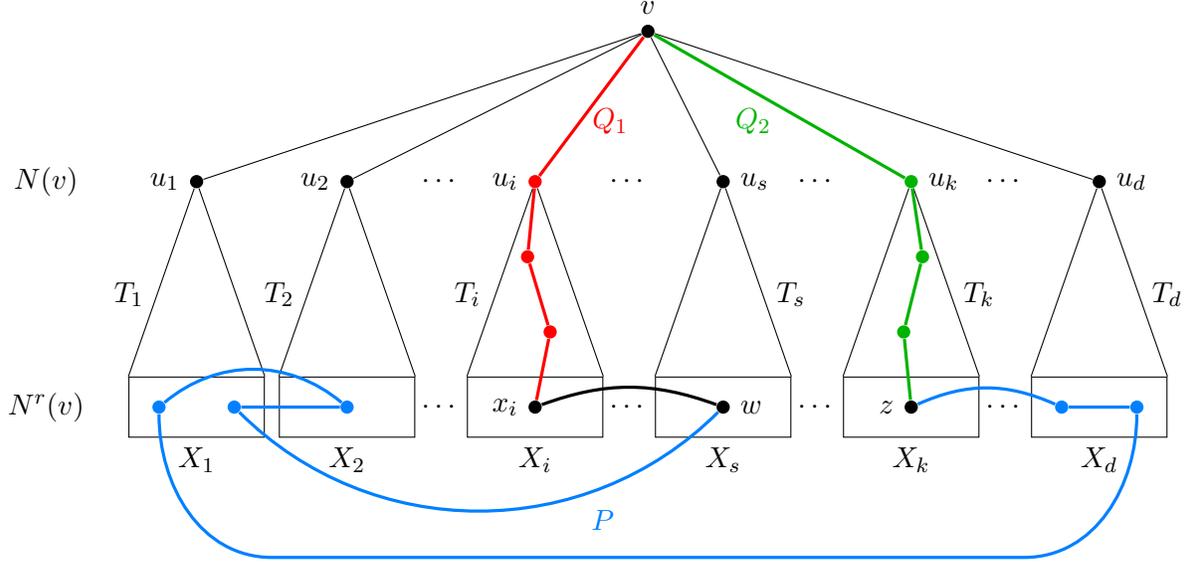
        
        Let $Q_1$ be the path from $v$ to $x_i$ in $T_i$, and $Q_2$ be the path from $v$ to $z$ in $T_k$. Both $Q_1$ and $Q_2$ have length $r$, and they intersect exactly on $\{v\}$. The concatenation $Q$ of $Q_1$ and $Q_2$ is thus a path of length $2r$ from $x_i$ to $z$, intersecting $H$ exactly on $\{x_i,z\}$. Hence, the concatenation of $P$ and $Q$ is a path of length exactly $2t$ from $w$ to $x_i$. Adding the edge $wx_i$ to this path yields a copy of $C_{2t+1}$ in $G$, a contradiction. The result follows.
    \end{proofclaim}

    For each $i\in [d]$ we let $H_i \coloneqq H[X_i]$. Observe that $H_i \subseteq L_{r-1}(u_i)$, and so by the induction hypothesis $H_i$ is $(2t-1)$-degenerate.
    For each $i<j$ we let $H_{ij}\coloneqq H[X_i,X_j]$ be the bipartite subgraph of $H$ induced by $X_i$ and $X_j$.
    
    \begin{claim}\label{claim:deg_Hij}
        For every pair of integers $i<j$, $H_{ij}$ is $(t-r)$-degenerate.
    \end{claim}
    \begin{proofclaim}
        Assume by contradiction that $H_{ij}$ contains a subgraph $H'$ with minimum degree $\delta(H') \ge t-r+1$. Thus, since $H'$ is bipartite, it contains a path $P_{ij}$ of length $2t-2r + 1$ with endpoints $y\in X_i$ and $z\in X_j$. Let $P_i$ (resp. $P_j$) be the path from $v$ to $y$ (resp. $z$) in $T$. Then $P_i \cup P_j \cup P_{ij}$ is a cycle of length $2t + 1$ in $G$, a contradiction. 
    \end{proofclaim}
     From Claim~\ref{claim:neighbourhood_Xj} we infer that every connected component in $H$ is contained in some $H_i$ or some $H_{ij}$. In either case, it is $(2t-1)$-degenerate, so it contains a vertex of degree less than $2t$, a contradiction.
\end{proof}

\subsection{Forbidding an even cycle}

We prove the following, which corresponds to the case of \Cref{lem:Lr(v)} where we forbid an even cycle.

\begin{lemma}
\label{lem:Lr(v)_even}
    Let $t\ge 1$ be a fixed integer, and let $G$ be a $C_{2t+2}$-free graph. Then for every vertex $v\in V(G)$ and every integer $1\le r\le t$, the subgraph $L_r(v)$ is $2t$-degenerate.
\end{lemma}

The proof of \Cref{lem:Lr(v)_even} relies on two intermediate results by Pikhurko in the study of the Tur\' an function of even cycles \cite{Pik12}. For an integer $k\ge 3$, a \emph{$\Theta_k$-graph} is a cycle of length at least $2k$ with a chord.

\begin{lemma}[Pikhurko {\cite[Lemma~2.2]{Pik12}}]
\label{lem:theta-degmin}
    Assume $k\ge 3$. Every bipartite graph of minimum degree at least $k$ contains a $\Theta_k$-subgraph.
\end{lemma}

\begin{lemma}[Pikhurko {\cite[Claim~3.1]{Pik12}}]
\label{lem:theta-free}
    Let $G$ be a $C_{2k}$-free graph, and let $v\in V(G)$. Then, for every integer $1\leq r\leq k-1$, $L_r(v)$ contains no bipartite $\Theta_k$-subgraph.
\end{lemma}

Using these two results, we are ready to prove \Cref{lem:Lr(v)_even}.

\begin{proof}[Proof of \Cref{lem:Lr(v)_even}]
    Assume for the sake of contradiction that there is a subgraph $H\subseteq L_r(v)$ of minimum degree $\delta(H) \ge 2t+1$, for some $1\le r \le t$.

    When $r=1$, this is a contradiction to Proposition~\ref{prop:C_ell_Nv_degenerate}.
    Henceforth, we assume that $r\ge 2$. A classical analysis of a maximum cut shows that every graph of minimum degree at least $2t+1$ contains a bipartite subgraph of minimum degree at least $t+1$. Applying Lemma~\ref{lem:theta-degmin} (with $k=t+1\geq 3$) to such a bipartite subgraph of $H$ yields a bipartite $\Theta_{t+1}$-subgraph. On the other hand, $G$ being $C_{2t+2}$-free, Lemma~\ref{lem:theta-free} ensures that $L_r(v)$ contains no bipartite $\Theta_{t+1}$-subgraph, a contradiction.
\end{proof}

\bibliography{discrepency}
\bibliographystyle{abbrv}

\end{document}